\documentclass[a4paper,11pt,reqno]{amsart}

\usepackage{dsfont,amssymb}

\numberwithin{equation}{section}

\newtheorem{theorem}{Theorem}[section]
\newtheorem{corollary}[theorem]{Corollary}
\newtheorem{lemma}[theorem]{Lemma}
\newtheorem{proposition}[theorem]{Proposition}
\theoremstyle{remark}

\newtheorem{example}[theorem]{Example}
\theoremstyle{definition}
\newtheorem{definition}[theorem]{Definition}

\newcommand\bp{\begin{proof}}
\newcommand\ep{\end{proof}}
\newcommand\ee{\nopagebreak\mbox{\ }\hfill$\diamondsuit$}

\newcommand\Dhat{{\hat\Delta}}

\newcommand\Ad{\operatorname{Ad}}
\newcommand\Aut{\operatorname{Aut}}
\newcommand\ev{\operatorname{ev}}

\newcommand{\C}{{\mathbb C}}
\newcommand{\R}{{\mathbb R}}
\newcommand\T{{\mathbb T}}
\newcommand\Z{{\mathbb Z}}

\newcommand{\A}{{\mathcal A}}
\newcommand{\F}{{\mathcal F}}
\newcommand{\M}{{\mathcal M}}

\newcommand\enu[1]{\smallskip\newline\makebox[5mm][l]{\rm(#1)}}

\begin{document}

\title{Deformation of operator algebras by Borel cocycles}

\author[J. Bhowmick]{Jyotishman Bhowmick}

\email{jyotishmanb@gmail.com}

\author[S. Neshveyev]{Sergey Neshveyev}

\email{sergeyn@math.uio.no}

\author[A. Sangha]{Amandip Sangha}

\email{amandip.s.sangha@gmail.com}

\address{Department of Mathematics, University of Oslo,
P.O. Box 1053 Blindern, NO-0316 Oslo, Norway}

\thanks{Supported by the Research Council of Norway.}

\date{July 11, 2012; minor corrections April 19, 2013}

\begin{abstract}
Assume that we are given a coaction $\delta$ of a locally compact group $G$ on a C$^*$-algebra~$A$ and a $\T$-valued Borel $2$-cocycle $\omega$ on $G$. Motivated by the approach of Kasprzak to Rieffel's deformation we define a deformation $A_\omega$ of $A$. Among other properties of $A_\omega$ we show that $A_\omega\otimes K(L^2(G))$ is canonically isomorphic to $A\rtimes_\delta\hat G\rtimes_{\hat\delta,\omega}G$. This, together with a slight extension of a result of Echterhoff et al., implies that for groups satisfying the Baum-Connes conjecture the $K$-theory of $A_\omega$ remains invariant under homotopies of $\omega$.
\end{abstract}

\maketitle


\section*{Introduction}

Assume $G$ is a discrete group and $A=\oplus_{g\in G}A_g$ is a $G$-graded algebra. Given a $\C^*$-valued $2$-cocycle~$\omega$ on $G$ we can define a new product on $A$ by the formula $a_g\cdot a_h=\omega(g,h)a_ga_h$ for $a_g\in A_g$ and $a_h\in A_h$. Some of the well-known examples of C$^*$-algebras, such as irrational rotation algebras and, more generally, twisted group C$^*$-algebras or twisted crossed products, are operator algebraic variants of this construction. Nevertheless the question what this construction means for a general C$^*$-algebra~$A$ and a locally compact group $G$ has no obvious answer. A natural replacement of a $G$-grading is a coaction of $G$ on $A$. But then the subspaces $A_g$ are often trivial for non-discrete $G$ and it is not clear how to define the new product.

In~\cite{Ri1} Rieffel succeeded in defining the product in the case $G=\R^d$ using oscillatory integrals. A few years ago Kasprzak~\cite{Kas} proposed an alternative approach that works for any locally compact group $G$ and a continuous $\T$-valued $2$-cocycle $\omega$. In fact, he considered only abelian groups and, correspondingly, actions of $\hat G$ rather than coactions of $G$, but it is easy to see that his construction makes sense for arbitrary $G$. It should also be mentioned that for discrete groups a different, but equivalent, approach has been recently suggested by Yamashita~\cite{Yam}.  Kasprzak's idea is as follows. Given a coaction $\delta$ of $G$ on $A$, consider the dual action $\hat\delta$ of $G$ on~$A\rtimes_\delta \hat G$. Using the cocycle $\omega$ we can deform this action to a new action $\hat\delta^\omega$. Then by general results on crossed products it turns out that $A\rtimes_\delta \hat G$ has another crossed product decomposition $A_\omega\rtimes_{\delta^\omega}\hat G$ such that~$\hat\delta^\omega$ becomes dual to $\delta^\omega$. The C$^*$-algebra $A_\omega$ is the deformation of $A$ we are looking for.

The goal of this note is to define $A_\omega$ for arbitrary Borel cocycles $\omega$. For abelian groups, restricting to continuous cocycles is not a serious omission, essentially since Borel cocycles correspond to Borel bicharacters and these are automatically continuous. But for general groups the class of continuous cocycles is too small and the right class is that of Borel cocycles~\cite{Moore1,Moore3}. Given a Borel cocycle $\omega$, there are no obvious reasons for the twisted dual action $\hat\delta^\omega$ to be well-defined on $A\rtimes_\delta \hat G$. What started this work is the observation that $\hat\delta^\omega$ is well-defined for dual coactions. Since any coaction is stably exterior equivalent to a dual coaction, and it is natural to expect that exterior equivalent coactions produce strongly Morita equivalent deformations, this suggested that $A_\omega$ could be defined for arbitrary~$\delta$. In the end, though, we found it easier to relate $A_\omega$ to twisted crossed products rather than to use dual coactions. This simplifies proofs, but the fundamental reasons for why $A_\omega$ is well-defined become somewhat hidden.

Our deformed algebras $A_\omega$ enjoy a number of expected properties. In particular, they come with canonical coactions $\delta^\omega$. However, the isomorphism $A\rtimes_\delta \hat G\cong A_\omega\rtimes_{\delta^\omega}\hat G$, which played an important role in~\cite{Kas} and~\cite{Yam}, is no longer available for general cocycles. Instead we construct an explicit isomorphism $A_\omega\otimes K(L^2(G))\cong A\rtimes_\delta\hat G\rtimes_{\hat\delta,\omega}G$, which is equally well suited for studying $A_\omega$.

Let us finally say a few words about sources of examples of coactions. The easiest is, of course, to take the dual coaction on a crossed product $A=B\rtimes_\alpha G$. In this case the deformation produces the twisted crossed product $B\rtimes_{\alpha,\omega}G$, as expected. But even if we start with dual coactions, we can get new coactions by taking e.g.\ free products. Given a corepresentation of the dual quantum group~$\hat G$, we can also consider infinite tensor products, as well as free Araki-Woods factors, see~\cite{V0} and references therein.

\medskip
\noindent
{\bf Acknowledgement.} It is our pleasure to thank Pawe{\l} Kasprzak and Makoto Yamashita for the inspiring correspondence.

\bigskip


\section{Actions, coactions and crossed products}

In this preliminary section we fix our notation and list a number of facts that we will freely use later.

Let $G$ be a second countable locally compact group. Fix a left-invariant Haar measure on $G$. Denote by $\lambda$ and $\rho$ the left and right regular representations on $G$. We will usually identify the reduced group C$^*$-algebra $C^*_r(G)$ with its image under~$\lambda$. Similarly, we will usually identify $C_0(G)$ with the algebra of operators of multiplication by functions on~$L^2(G)$. Denote by $K$ the algebra of compact operators on $L^2(G)$.

Denote by $\Delta\colon C_0(G)\to M(C_0(G)\otimes C_0(G))=C_b(G\times G)$ and $\Dhat\colon C_r^*(G)\to M(C^*_r(G)\otimes C^*_r(G))$ the standard comultiplications, so
$$
\Delta(f)(g,h)=f(gh),\ \ \Dhat(\lambda_g)=\lambda_g\otimes\lambda_g.
$$

Let $W\in M(C_0(G)\otimes C^*_r(G))$ be the fundamental unitary, defined by
$$
(W\xi)(s,t)=\xi(s,s^{-1}t)\ \ \text{for}\ \ \xi\in L^2(G\times G).
$$
In other words, if we identify $M(C_0(G)\otimes C^*_r(G))$ with the algebra of bounded strictly continuous maps $G\to M(C_r^*(G))$, then $W(g)=\lambda_g$. We have
$$
W^*(1\otimes f)W=\Delta(f)\ \ \text{for}\ f\in C_0(G),
$$
$$
W(\lambda_g\otimes1)W^*=\Dhat(\lambda_g)\ \ \text{and}\ \
W^*(\rho_g\otimes1)W=\rho_g\otimes\lambda_g\ \text{for}\ \ g\in G.
$$
We will also use the unitary $V=(\rho\otimes\iota)(W_{21})\in M(\rho(C^*_r(G))\otimes C_0(G))$. We have
$$
V(f\otimes 1)V^*=\Delta(f)\ \ \text{for}\ f\in C_0(G),
$$
$$
V^*(1\otimes\rho_g)V=\rho_g\otimes\rho_g\ \ \text{and}\ \
V(1\otimes\lambda_g)V^*=\rho_g\otimes\lambda_g\ \text{for}\ \ g\in G.
$$

\medskip

Assume now that $\alpha\colon G\to\Aut(B)$ is a (continuous) action of $G$ on a C$^*$-algebra $B$. We consider~$\alpha$ as a homomorphism $\alpha\colon B\to M(B\otimes C_0(G))$, so that $\alpha(b)(g)=\alpha_g(b)$. Then $(\alpha\otimes\iota)\alpha=(\iota\otimes\Delta)\alpha$. We define the reduced crossed product $B\rtimes_\alpha G$ by
$$
B\rtimes_\alpha G=\overline{\alpha(B)(1\otimes\rho(C^*_r(G)))}\subset M(B\otimes K).
$$
This is equivalent to the standard definition. Since {\bf we consider only reduced crossed products in this paper}, we omit $r$ in the notation.

By a coaction of $G$ on a C$^*$-algebra $A$ we mean a non-degenerate injective homomorphism $\delta\colon A\to M(A\otimes C^*_r(G))$ such that $(\delta\otimes\iota)\delta=(\iota\otimes\Dhat)\delta$ and the space $\delta(A)(1\otimes C^*_r(G))$ is a dense subspace of~$A\otimes C^*_r(G)$. The crossed product is then defined by
$$
A\rtimes_\delta\hat G=\overline{\delta(A)(1\otimes C_0(G))}\subset M(A\otimes K).
$$
The algebra $A\rtimes_\delta\hat G$ is equipped with the dual action $\hat\delta$ of $G$ defined by $\hat\delta_g=\Ad(1\otimes \rho_g)$. Thinking of $\hat\delta$ as a homomorphism $A\rtimes_\delta\hat G\to M((A\rtimes_\delta\hat G)\otimes C_0(G))$, we have
$$
\hat\delta(\delta(a))=\delta(a)\otimes 1, \ \ \hat\delta(1\otimes f)=1\otimes\Delta(f).
$$
It follows that
$$
\hat\delta(x)=V_{23}(x\otimes 1)V_{23}^*\ \ \text{for}\ \ x\in A\rtimes_\delta\hat G\subset M(A\otimes K).
$$

Similarly, starting with an action $\alpha$ of $G$ on $B$ we get a dual coaction $\hat\alpha$ of $G$ on $B\rtimes_\alpha G$ such that
$$
\hat\alpha(\alpha(b))=\alpha(b)\otimes1,\ \ \hat\alpha(1\otimes\rho_g)=1\otimes\rho_g\otimes\lambda_g.
$$
Therefore
$$
\hat\alpha(x)=W^*_{23}(x\otimes 1)W_{23}\ \ \text{for}\ \ x\in B\rtimes_\alpha G\subset M(B\otimes K).
$$

\medskip

A $1$-cocycle for an action $\alpha$ of $G$ on $B$ is a strictly continuous family $U=\{u_g\}_{g\in G}$ of unitaries in~$M(B)$ such that $u_{gh}=u_g\alpha_g(u_h)$. Given such a cocycle, we can define a new action $\alpha_U$ of~$G$ on~$B$ by $\alpha_{U,g}=u_g\alpha_g(\cdot)u_g^*$. The actions $\alpha$ and $\alpha_U$ are called exterior equivalent. We have an isomorphism $B\rtimes_\alpha G\cong B\rtimes_{\alpha_U} G$ respecting the dual coactions, defined~by
$$
\alpha(b)\mapsto\alpha_U(b),\ \ 1\otimes \rho_g\mapsto\alpha_U(u_g^*)(1\otimes\rho_g).
$$
If we think of $U$ as an element of $M(B\otimes C_0(G))$, then this isomorphism is implemented by the inner automorphism $\Ad U$ of $M(B\otimes K)$.

Similarly, a $1$-cocycle for a coaction $\delta$ of $G$ on $A$ is a unitary $U\in M(A\otimes C^*_r(G))$ such that $(\iota\otimes\Dhat)(U)=(U\otimes1)(\delta\otimes\iota)(U)$. Given such a cocycle, we can define a new coaction $\delta_U$ by $\delta_U(a)=U\delta(a)U^*$. The coactions $\delta$ and $\delta_U$ are called exterior equivalent. The inner automorphism~$\Ad U$ of~$M(A\otimes K)$ defines an isomorphism of $A\rtimes_\delta\hat G $ onto $A\rtimes_{\delta_U}\hat G$ respecting the dual actions, see~\cite[Theorem~2.9]{LPRS}.

In particular, given a coaction $\delta$ of $G$ on $A$ we can consider the coaction $a\otimes T\mapsto \delta(a)_{13}(1\otimes T\otimes1)$ of $G$ on $A\otimes K$, then take the $1$-cocycle $1\otimes W^*$ for this coaction (the cocycle identity means that $(\iota\otimes\Dhat)(W)=W_{13}W_{12}$) and get a new coaction on $A\otimes K$. In order to lighten the notation we will denote this new coaction by $\delta_{W^*}$. Then the Takesaki-Takai(-Katayama-Baaj-Skandalis) duality states that
$$
(A\rtimes_\delta \hat G\rtimes_{\hat\delta}G,\hat{\hat\delta})\cong (A\otimes K,\delta_{W^*}).
$$
Explicitly, the isomorphism is given by
$$
\hat\delta(\delta(a))=\delta(a)\otimes1\mapsto\delta(a), \ \ \hat\delta(1\otimes f)=1\otimes\Delta(f)\mapsto 1\otimes f, \ \ 1\otimes1\otimes\rho_g\mapsto 1\otimes\rho_g.
$$
If we identify $A\otimes K$ with $\delta(A)\otimes K\subset M(A\otimes K\otimes K)$, then this isomorphism is simply $\Ad W_{23}$.

\medskip

We finish this section by discussing how to recover $A$ from $A\rtimes_\delta\hat G$ for a coaction $\delta$. Consider the homomorphism
$$
\eta\colon A\rtimes_\delta \hat G\to M((A\rtimes_\delta\hat G)\otimes K)\subset M(A\otimes K\otimes K)
$$
defined by $\eta(x)=W_{23}\hat\delta(x)W_{23}^*$. In other words, $\eta$ is the composition of $\hat\delta\colon A\rtimes_\delta \hat G\to M(A\rtimes_\delta \hat G\rtimes_{\hat\delta}G)$ with the Takesaki-Takai duality isomorphism $A\rtimes_\delta \hat G\rtimes_{\hat\delta}G\cong \delta(A)\otimes K$. Explicitly,
$$
\eta(\delta(a))=(\delta\otimes\iota)\delta(a),\ \ \eta(1\otimes f)=1\otimes f\in M((A\rtimes_\delta\hat G)\otimes K)
$$
From this we see that $\delta(A)\subset M(A\rtimes_\delta \hat G)$ is the closed linear span of elements of the form $(\iota\otimes\varphi)\eta(x)$ with $x\in A\rtimes_\delta\hat G$ and $\varphi\in K^*$.

More generally, assume we are given an action $\alpha$ of $G$ on a C$^*$-algebra $B$ and a nondegenerate homomorphism $\pi\colon C_0(G)\to M(B)$ such that $\alpha(\pi(f))=(\pi\otimes\iota)\Delta(f)$. Put $X=(\pi\otimes\iota)(W)$ and consider the homomorphism
$$
\eta\colon B\to M(B\otimes K), \ \ \eta(x)=X\alpha(x)X^*.
$$
Then by a Landstad-type result of Quigg \cite[Theorem~3.3]{Qu} and, more generally, Vaes~\cite[Theorem~6.7]{V}, the closed linear span $A\subset M(B)$ of elements of the form  $(\iota\otimes\varphi)\eta(x)$, with $x\in B$ and $\varphi\in K^*$, is a C$^*$-algebra, the formula $\delta(a)=X(a\otimes1)X^*$ defines a coaction of $G$ on $A$, and $\eta$ becomes an isomorphism $B\cong A\rtimes_\delta \hat G$ that intertwines $\alpha$ with $\hat\delta$.


\bigskip

\section{Deformation of algebras}

Denote by $Z^2(G;\T)$ the set of $\T$-valued Borel $2$-cocycles on $G$, so $\omega\in Z^2(G;\T)$ is a Borel function $G\times G\to \T$ such that
$$
\omega(g,h)\omega(gh,k)=\omega(g,hk)\omega(h,k).
$$
For every cocycle $\omega$ consider also the cocycles $\tilde\omega$ and $\bar\omega$ defined by $$\tilde\omega(g,h)=\omega(h^{-1},g^{-1})\ \ \text{and}\ \ \bar\omega(g,h)=\overline{\omega(g,h)}.$$

Define operators $\lambda^\omega_g$ and $\rho^{\tilde\omega}_g$ on~$L^2(G)$~by
\footnote{The operators $\lambda^\omega_g$ and $\rho^{\tilde\omega}_g$ are more commonly defined by $\lambda^\omega_g=\lambda_g\omega(g,\cdot)=\omega(g,g^{-1}\cdot)\lambda_g$ and $\rho^{\tilde\omega}_g=\rho_g\omega(\cdot,g^{-1})=\omega(\cdot g,g^{-1})\rho_g$.
With our definition some of the formulas will look better.
If the cocycle $\omega$ satisfies $\omega(g,e)=\omega(e,g)=\omega(g,g^{-1})=1$ for all $g\in G$, then the two definitions coincide, that is to say $\omega(h^{-1},g)=\omega(g,g^{-1}h)$, which follows by applying the cocycle identity for $\omega$ to the triple $(h^{-1},g,g^{-1}h)$. Any cocycle is cohomologous to a cocycle satisfying the above normalization conditions, so in principle we could consider only such cocycles.}
$$
\lambda^\omega_g=\tilde\omega(g^{-1},\cdot)\lambda_g,\ \ \rho^{\tilde\omega}_g=\tilde\omega(\cdot,g)\rho_g.
$$
Then
$$
\lambda^\omega_g\lambda^\omega_h=\omega(g,h)\lambda^\omega_{gh},\ \
\rho^{\tilde\omega}_g\rho^{\tilde\omega}_h=\tilde\omega(g,h)\rho^{\tilde\omega}_{gh}\ \ \text{and}\ \
[\lambda^\omega_g,\rho^{\tilde\omega}_h]=0\ \ \text{for all}\ \ g,h\in G.
$$

\medskip

Fix now a cocycle $\omega\in Z^2(G;\T)$ and consider a coaction $\delta$ of $G$ on a C$^*$-algebra $A$. Assume first that the cocycle $\omega$ is continuous. In this case the functions $\tilde\omega(\cdot,g)$ belong to the multiplier algebra of $C_0(G)$, so we can define a new twisted dual action $\hat\delta^\omega$ on $A\rtimes_\delta\hat G$ by letting $\hat\delta^\omega_g=\Ad(1\otimes\rho^{\tilde\omega}_g)$. In other words, if we consider $\tilde\omega$ as a multiplier of $C_0(G)\otimes C_0(G)$, then
$$
\hat\delta^\omega(x)=\tilde\omega_{23}\hat\delta(x){\tilde\omega}^*_{23}
=\tilde\omega_{23}V_{23}(x\otimes1)V_{23}^*{\tilde\omega}^*_{23}
\in M(A\otimes K\otimes K).
$$
For $f\in C_0(G)$ we obviously have $\hat\delta^\omega(1\otimes f)=\hat\delta(1\otimes f)=1\otimes\Delta(f)$. By the Landstad-type duality result of Quigg and Vaes it follows that $\hat\delta^\omega$ is the dual action on a crossed product $A_\omega\rtimes_{\delta^\omega}\hat G$ for some C$^*$-subalgebra $A_\omega\subset M(A\rtimes_\delta\hat G)\subset M(A\otimes K)$ and a coaction $\delta^\omega$ of $G$, and this subalgebra is defined using slice maps applied to the image of $A\rtimes_\delta\hat G$ under the homomorphism
$$
\eta^\omega\colon A\rtimes_\delta \hat G\to M(A\otimes K\otimes K),\ \
\eta^\omega(x)=W_{23}\tilde\omega_{23}\hat\delta(x){\tilde\omega}^*_{23}W_{23}^*.
$$

If the cocycle $\omega$ is only assumed to be Borel, it is not clear whether the action $\hat\delta^\omega$ is well-defined. Nevertheless, the homomorphism $\eta^\omega\colon A\rtimes_\delta G\to M(A\otimes K\otimes K)$ defined above still makes sense. Therefore we can give the following definition.

\begin{definition}
The $\omega$-deformation of a C$^*$-algebra $A$ equipped with a coaction $\delta$ of $G$ is the C$^*$-subalgebra $A_\omega\subset M(A\otimes K)$ generated by all elements of the form $$(\iota\otimes\iota\otimes\varphi)\eta^\omega\delta(a)
=(\iota\otimes\iota\otimes\varphi)\Ad(W_{23}\tilde\omega_{23})(\delta(a)\otimes1),$$ where $a\in A$ and $\varphi\in K^*$.
\end{definition}

In case we want to stress that the deformation is defined using the coaction $\delta$, we will write $A_{\delta,\omega}$ instead of $A_\omega$.

Note that if we considered elements of the form $(\iota\otimes\iota\otimes\varphi)\eta^\omega(x)$ for all $x\in A\rtimes_\delta\hat G$, this would not change the algebra $A_\omega$, since $\eta^\omega(1\otimes f)=1\otimes1\otimes f$.

\smallskip

In order to get an idea about the structure of $A_\omega$ consider the C$^*$-algebra  $C^*_r(G,\omega)$ generated by operators of the form
$$
\lambda^\omega_f=\int_Gf(g)\lambda^\omega_gdg, \ \ f\in L^1(G).
$$
When necessary we denote by $\lambda^\omega$ the identity representation of $C^*_r(G,\omega)$ on $L^2(G)$. A simple computation shows that
\begin{equation} \label{ebasic}
W\tilde\omega(\lambda_g\otimes1)\tilde\omega^*W^*=\lambda^\omega_g\otimes \lambda^{\bar\omega}_g.
\end{equation}
The map $g\mapsto \lambda^\omega_g\otimes \lambda^{\bar\omega}_g$ therefore defines a representation of $G$ on $L^2(G\times G)$ that is quasi-equivalent to the regular representation, so it defines a representation of $C^*_r(G)$. Denote this representation by~$\lambda^\omega\boxtimes\lambda^{\bar\omega}$. We can then write
$$
\eta^\omega\delta(a)=(\iota\otimes(\lambda^\omega\boxtimes\lambda^{\bar\omega}))\delta(a)\ \ \text{for}\ \ a\in A.
$$
Since the image of $C^*_r(G)$ under $\lambda^\omega\boxtimes\lambda^{\bar\omega}$ is contained in $M(C^*_r(G,\omega)\otimes C^*_r(G,\bar\omega))$, we see in particular that $A_\omega\subset M(A\otimes C^*_r(G,\omega))$.

\begin{example} \label{exdiscr}
Assume the group $G$ is discrete. Denote by $A_g\subset A$ the spectral subspace corresponding to $g\in G$, so $A_g$ consists of all elements $a\in A$ such that $\delta(a)=a\otimes\lambda_g$. The spaces $A_g$, $g\in G$, span a dense $*$-subalgebra $\A\subset A$.
By \eqref{ebasic}, if $a\in A_g$ then
$
\eta^\omega\delta(a)=a\otimes\lambda^\omega_g\otimes \lambda^{\bar\omega}_g.
$
Thus the linear span of elements $(\iota\otimes\iota\otimes\varphi)\eta^\omega\delta(a)$, with $a\in\A$ and $\varphi\in K^*$, coincides with the linear span $\A_\omega$ of elements $a\otimes\lambda^\omega_g$, with $a\in A_g$ and $g\in G$. The space $\A_\omega$ is already a $*$-algebra and $A_\omega$ is the closure of $\A_\omega$ in $A\otimes C^*_r(G,\omega)$. In particular, we see that for discrete groups our definition of $\omega$-deformation is equivalent to that of Yamashita, see~\cite[Proposition~2]{Yam}.
\ee
\end{example}

The following theorem is the first principal result of this section.

\begin{theorem} \label{tmain}
The C$^*$-algebra $A_\omega\subset M(A\otimes K)$ coincides with the norm closure of the linear span of elements of the form $(\iota\otimes\iota\otimes\varphi)\eta^\omega\delta(a)$, where $a\in A$ and $\varphi\in K^*$.
\end{theorem}

While proving this theorem we will simultaneously obtain a description of $A_\omega\otimes K$. We need to introduce more notation in order to formulate the result.

\smallskip

In addition to $\lambda^\omega$ we have another equivalent representation $\rho^\omega$ of $C^*_r(G,\omega)$ on $L^2(G)$ that maps $\lambda^\omega_g\in M(C^*_r(G,\omega))$ into $\rho^\omega_g$.

Given an action $\alpha$ of $G$ on a C$^*$-algebra $B$, the reduced twisted crossed product is defined by
$$
B\rtimes_{\alpha,\omega}G=\overline{\alpha(B)(1\otimes\rho^\omega(C^*_r(G,\omega)))}\subset M(B\otimes K).
$$
The reduced twisted crossed product has a dual coaction, which we again denote by $\hat\alpha$, defined by
$$
\hat\alpha(x)=W_{23}^*(x\otimes1)W_{23},\ \ \text{so}\ \ \hat\alpha(\alpha(b))=\alpha(b)\otimes1, \ \ \hat\alpha(1\otimes\rho^\omega_g)=1\otimes\rho^\omega_g\otimes\lambda_g.
$$

The last ingredient that we need is the well-known fact that the cocycles $\tilde\omega$ and $\bar\omega$ are cohomologous. Explicitly,
$$
\tilde\omega(g,h)=\bar\omega(g,h)v(g)v(h)v(gh)^{-1},\ \ \text{where}\ \ v(g)=\omega(g^{-1},g)\omega(e,e).
$$
This follows from the cocycle identities
$$
\omega(h^{-1},g^{-1})\omega(h^{-1}g^{-1},gh)=\omega(h^{-1},h)\omega(g^{-1},gh), \ \
\omega(g^{-1},gh)\omega(g,h)=\omega(g^{-1},g)\omega(e,h);
$$
recall also that $\omega(e,h)=\omega(e,e)$ for all $h$, which follows from the cocycle identity applied to $(e,e,h)$.

\medskip

We can now formulate our second principal result.

\begin{theorem} \label{tmain2}
Put $u(g)=\overline{\omega(g^{-1},g)\omega(e,e)}$. Then the map
$$
\Ad((1\otimes W\tilde\omega)(1\otimes1\otimes u))\colon A\rtimes_\delta\hat G\rtimes_{\hat\delta,\omega}G\to M(A\otimes K\otimes K)
$$
defines an isomorphism $A\rtimes_\delta\hat G\rtimes_{\hat\delta,\omega}G\cong A_{\omega}\otimes K$.
\end{theorem}

For discrete groups the fact that the C$^*$-algebras $A_\omega$ and $A\rtimes_\delta\hat G\rtimes_{\hat\delta,\omega}G$ are strongly Morita equivalent was observed by Yamashita~\cite[Corollary~15]{Yam}.

\bp[Proof of Theorems~\ref{tmain} and~\ref{tmain2}]
Denote by $\theta$ the map in the formulation of Theorem~\ref{tmain2}. In order to compute its image, observe first that since $\bar{\tilde\omega}(h,g)=\omega(h,g)u(h)u(g)u(hg)^{-1}$, we have
$$
u\rho^\omega_gu^*=\overline{u(g)}\rho^{\bar{\tilde\omega}}_g.
$$
Next, it is straightforward to check that $W\tilde\omega$ commutes with $1\otimes\rho^{\bar{\tilde\omega}}_g$. We thus see that $\theta$ acts as
$$
\delta(a)\otimes1\mapsto\eta^{\omega}\delta(a),\ \
1\otimes\Delta(f)\mapsto 1\otimes1\otimes f, \ \ 1\otimes1\otimes\rho^\omega_g\mapsto 1\otimes1\otimes u\rho^\omega_gu^*.
$$
In particular, we see that the image of the C$^*$-subalgebra
$$
\overline{(1\otimes\Delta(C_0(G)))(1\otimes1\otimes\rho^\omega(C^*_r(G,\omega)))}\cong C_0(G)\rtimes_{\Ad\rho,\omega}G
$$
of $M(A\rtimes_\delta\hat G\rtimes_{\hat\delta,\omega}G)$ is
$$
1\otimes1\otimes\overline{uC_0(G)C^*_r(G,\omega)u^*}=1\otimes1\otimes K.
$$
Therefore $1\otimes1\otimes K$ is a nondegenerate C$^*$-subalgebra of $M(\theta(A\rtimes_\delta\hat G\rtimes_{\hat\delta,\omega}G))\subset M(A\otimes K\otimes K)$.
It follows that there exists a uniquely defined C$^*$-subalgebra $A_1\subset M(A\otimes K)$ such that
$$
\theta(A\rtimes_\delta\hat G\rtimes_{\hat\delta,\omega}G)=A_1\otimes K.
$$
By definition of crossed products and the above computation of $\theta$ we then have
$$
A_1\otimes K=\overline{\eta^\omega\delta(A)(1\otimes1\otimes K)}.
$$
Applying the slice maps $\iota\otimes\iota\otimes\varphi$ we conclude that the closed linear span of elements of the form $(\iota\otimes\iota\otimes\varphi)\eta^\omega\delta(a)$ coincides with the C$^*$-algebra $A_1$. This finishes the proof of both theorems.
\ep

Theorem~\ref{tmain2} essentially reduces the study of $\omega$-deformations to that of (twisted) crossed products. As a simple illustration let us prove the following result that refines and generalizes \cite[Proposition~14]{Yam}.

\begin{proposition}
Assume we are given two exterior equivalent coactions $\delta$ and $\delta_X$ of $G$ on a C$^*$-algebra $A$. Then $A_{\delta,\omega}\otimes K\cong A_{\delta_X,\omega}\otimes K$.
\end{proposition}

\bp Since $\delta$ and $\delta_X$ are exterior equivalent, we have $(A\rtimes_{\delta}\hat G,\hat\delta)\cong(A\rtimes_{\delta_X}\hat G,\hat\delta_X)$, and hence $A\rtimes_{\delta}\hat G\rtimes_{\hat\delta,\omega}G\cong A\rtimes_{\delta_X}\hat G\rtimes_{\hat\delta_X,\omega}G$.
\ep

Note that for continuous cocycles this result is also a consequence of the following useful fact combined with the Takesaki-Takai duality.

\begin{proposition}
If the cocycle $\omega$ is continuous, then any two exterior equivalent coactions have exterior equivalent twisted dual actions. More precisely, assume $X\in M(A\otimes C_r^*(G))$ is a $1$-cocycle for a coaction $\delta$ of $G$ on~$A$. Then the element
$U=X_{12}\tilde\omega_{23}X_{12}^*{\tilde\omega}^*_{23}\in M(A\otimes K\otimes C_0(G))$ is a $1$-cocycle for the action $\hat\delta^\omega_X$ of~$G$ on~$A\rtimes_{\delta_X}\hat G$, and the isomorphism $\Ad X\colon A\rtimes_\delta\hat G\to A\rtimes_{\delta_X}\hat G$ intertwines $\hat\delta^\omega$ with~$(\hat\delta^\omega_X)_U$.
\end{proposition}

\bp Denote by $\Psi$ the isomorphism $\Ad X\colon A\rtimes_\delta\hat G\to A\rtimes_{\delta_X}\hat G$ and put $$Y=1\otimes\tilde\omega\in M(1\otimes C_0(G)\otimes C_0(G))\subset M((A\rtimes_\delta\hat G)\otimes C_0(G))\cap M((A\rtimes_{\delta_X}\hat G)\otimes C_0(G)).$$
Then $U=(\Psi\otimes\iota)(Y)Y^*\in M((A\rtimes_{\delta_X}\hat G)\otimes C_0(G))$. In order to show that $U$ is a $1$-cocycle for $\hat\delta^\omega_X$, observe first that
\begin{equation} \label{e1}
(Y\otimes1)(\hat\delta_X\otimes\iota)(Y)
=(\iota\otimes\iota\otimes\Delta)(Y)\tilde\omega_{34},
\end{equation}
which is simply the cocycle identity for $\tilde\omega$. We also have the same identity for $\hat\delta$. Furthermore, since~$\Psi$ intertwines $\hat\delta$ with $\hat\delta_X$, we also get
$$
((\Psi\otimes\iota )(Y)\otimes1)
(\hat\delta_X\otimes\iota)(\Psi\otimes\iota)(Y)
=(\iota\otimes\iota\otimes\Delta)(\Psi\otimes\iota)(Y)\tilde\omega_{34}.
$$
Multiplying this identity by the adjoint of \eqref{e1} we obtain
$$
((\Psi\otimes\iota )(Y)\otimes1)(\hat\delta_X\otimes\iota)(U)(Y^*\otimes 1)
=(\iota\otimes\iota\otimes\Delta)(U).
$$
Since $\hat\delta^\omega_X=Y\hat\delta_X(\cdot)Y^*$, this is exactly the cocycle identity
$$
(U\otimes1)(\hat\delta_X\otimes\iota)(U)=(\iota\otimes\iota\otimes\Delta)(U).
$$

Since $\hat\delta^\omega=Y\hat\delta(\cdot)Y^*$, $\hat\delta^\omega_X=Y\hat\delta_X(\cdot)Y^*$ and $\Psi$ intertwines $\hat\delta$ with $\hat\delta_X$, we immediately see that $\Psi$ intertwines $\hat\delta^\omega$ with $(\Psi\otimes\iota)(Y)\hat\delta_X(\cdot)(\Psi\otimes\iota)(Y)^*=U\hat\delta^\omega_X(\cdot)U^*$.
\ep

We finish the section with the following simple observation.

\begin{proposition}
Assume $\omega_1,\omega_2\in Z^2(G;\T)$ are cohomologous cocycles. Then $A_{\omega_1}\cong A_{\omega_2}$.
\end{proposition}

\bp By assumption there exists a Borel function $v\colon G\to\T$ such that $$\tilde\omega_1(g,h)=\tilde\omega_2(g,h)v(g)v(h)v(gh)^{-1},$$ that is, $\tilde\omega_1=\tilde\omega_2(v\otimes v)\Delta(v)^*$. Note that then $\lambda^{\omega_1}_g=v(g^{-1})v\lambda^{\omega_2}_gv^*$. Using that $W\Delta(v)W^*=1\otimes v$ and that $W$ commutes with $v\otimes1$, for any operator $x$ on $L^2(G)$ we get
$$
W\tilde\omega_1(x\otimes1)\tilde\omega_1^*W^*=(v\otimes v^*)W\tilde\omega_2(x\otimes1)\tilde\omega_2^*W^*(v^*\otimes v).
$$
This shows that
$$
\eta^{\omega_1}=\Ad(1\otimes v\otimes v^*)\eta^{\omega_2},
$$
which in turn gives $A_{\omega_1}=\Ad(1\otimes v)(A_{\omega_2})$.
\ep


\bigskip

\section{Canonical and dual coactions}

By the Landstad-type result of Quigg and Vaes the twisted dual action~$\hat\delta^\omega$, when it is defined, is dual to some coaction. The action $\hat\delta^\omega$ is apparently not always well-defined on $A\rtimes_\delta\hat G$. Nevertheless the new coaction on $A_\omega$ always makes sense.

\begin{theorem}
For any cocycle $\omega\in Z^2(G;\T)$ and a coaction $\delta$ of $G$ on a C$^*$-algebra $A$ we have:
\enu{i} the formula $\delta^\omega(x)=W_{23}(x\otimes1)W_{23}^*$ defines a coaction of~$G$ on $A_\omega$;
\enu{ii} if the twisted dual action $\hat\delta^\omega$ is well-defined on $A\rtimes_\delta G$, then $A\rtimes_\delta \hat G=\overline{A_\omega(1\otimes C_0(G))}$ and the map $\eta^\omega\colon A\rtimes_\delta \hat G\to M(A\otimes K\otimes K)$ gives an isomorphism $A\rtimes_\delta \hat G\cong A_\omega\rtimes_{\delta^\omega} \hat G$ that intertwines the twisted dual action~$\hat\delta^\omega$ on $A\rtimes_\delta \hat G$ with the dual action to~$\delta^\omega$ on~$A_\omega\rtimes_{\delta^\omega} \hat G$.
\end{theorem}

\bp (i) We repeat the computations of Vaes in the proof \cite[Theorem~6.7]{V}. Since
$$
W_{13}W_{12}=(\iota\otimes\Dhat)(W)=W_{23}W_{12}W_{23}^*,
$$
for $x=(\iota\otimes\iota\otimes\varphi)\eta^\omega(y)$, $y\in A\rtimes_\delta\hat G$, we have
\begin{align*}
\delta^\omega(x)&=(\iota\otimes\iota\otimes\varphi\otimes\iota)(W_{24}(\eta^\omega(y)\otimes1)W_{24}^*)\\
&=(\iota\otimes\iota\otimes\varphi\otimes\iota)(W_{24}W_{23}\tilde\omega_{23}(\hat\delta(y)\otimes1)
{\tilde\omega}^*_{23}W_{23}^*W_{24}^*)\\
&=(\iota\otimes\iota\otimes\varphi\otimes\iota)(W_{34}W_{23}W_{34}^*\tilde\omega_{23}(\hat\delta(y)\otimes1)
{\tilde\omega}^*_{23}W_{34}W_{23}^*W_{34}^*)\\
&=(\iota\otimes\iota\otimes\varphi\otimes\iota)(W_{34}(\eta^\omega(y)\otimes1)W_{34}^*).
\end{align*}
From this one can easily see that the closure of $\delta^\omega(A_\omega)(1\otimes1\otimes C^*_r(G))$ coincides with $A_\omega\otimes C^*_r(G)$, because $\overline{(K\otimes1)W(1\otimes C^*_r(G))}=K\otimes C^*_r(G)$ and $W^*(K\otimes C_r^*(G))=K\otimes C_r^*(G)$. Since $1\otimes W$ is a $1$-cocycle for the trivial coaction on $A\otimes K$ (so $(\iota\otimes\Dhat)(W)=W_{12}W_{13}$), the identity $(\iota\otimes\Dhat)\delta^\omega=(\delta^\omega\otimes\iota)\delta^\omega$ follows.

\smallskip

(ii) This is \cite[Theorem~6.7]{V} applied to the action $\hat\delta^\omega$.
\ep

The twisted dual action is well-defined for continuous cocycles, but as the following result shows it can also be well-defined even if the cocycle is only Borel.

\begin{proposition}
If $\delta$ is a dual coaction, then the twisted dual action $\hat\delta^\omega$ of $G$ on $A\rtimes_\delta\hat G$ is well-defined for any $\omega\in Z^2(G;\T)$.
\end{proposition}

\bp
By assumption we have $A=B\rtimes_\alpha G$ and $\delta=\hat\alpha$ for some $B$ and $\alpha$. Then $A\rtimes_\delta\hat G=B\rtimes_\alpha G\rtimes_{\hat\alpha}\hat G$ is the closure of
$$
(\alpha(B)\otimes1)(1\otimes(\rho\otimes\lambda)\Dhat(C^*_r(G)))(1\otimes1\otimes C_0(G))\subset M(B\otimes K\otimes K).
$$
We have to check that the inner automorphisms $\Ad(1\otimes1\otimes\rho^{\tilde\omega}_g)$ of $B\otimes K\otimes K$ define a (continuous) action of $G$ on this closure. Since these automorphisms act trivially on $\alpha(B)\otimes1$, we just have to check that the automorphisms $\Ad(1\otimes\rho^{\tilde\omega}_g)$ of $K\otimes K$ define an action on the C$^*$-algebra $$\overline{(\rho\otimes\lambda)\Dhat(C^*_r(G))(1\otimes C_0(G))}\cong C^*_r(G)\rtimes\hat G.$$
The operator $V$ commutes with $1\otimes\tilde\omega(\cdot,g)$, and $\Ad V^*$ maps the above algebra onto $1\otimes K$. Hence $\Ad(1\otimes\tilde\omega(\cdot,g))$, and therefore also $\Ad(1\otimes\rho^{\tilde\omega}_g)$, is a well-defined automorphism of that algebra. Finally, the continuity of the action holds, since any Borel homomorphism of $G$ into a Polish group, such as the group $\Aut(K)$, is automatically continuous.
\ep

For dual coactions it is, however, straightforward to describe the deformed algebra, see \cite[Example~8]{Yam} for the discrete group case. In order to formulate the result, define a unitary $W^\omega$ on $L^2(G\times G)$~by
$$
(W^\omega\xi)(g,h)=\tilde\omega(g^{-1},h)\xi(g,g^{-1}h).
$$
In other words, if we let $W^*(G,\omega)=C^*_r(G,\omega)''$, then $W^\omega\in L^\infty(G)\bar\otimes W^*(G,\omega)=L^\infty(G;W^*(G,\omega))$ and $W^\omega(g)=\lambda^\omega_g$.

\begin{proposition}
Assume $\alpha$ is an action of $G$ on a C$^*$-algebra $B$. Consider the dual coaction $\delta$ on $A=B\rtimes_\alpha G$. Then for any $\omega\in Z^2(G;\T)$ the map
$$
B\rtimes_{\alpha,\omega}G\mapsto M(B\otimes K\otimes K), \ \ x\mapsto W^{\omega*}_{23}(x\otimes1)W^\omega_{23},
$$
defines an isomorphism $(B\rtimes_{\alpha,\omega}G,\hat\alpha)\cong(A_\omega,\delta^\omega)$.
\end{proposition}

\bp First of all observe that by \eqref{ebasic} we have
$$
\eta^\omega(\delta(1\otimes\rho_g))=1\otimes \rho_g\otimes \lambda^\omega_g\otimes\lambda^{\bar\omega}_g.
$$
This implies that $A_\omega$ is the closed linear span of elements of the form
$$
(\delta(b)\otimes 1)\int_Gf(g)(1\otimes\rho_g\otimes \lambda^\omega_g)dg,
$$
where $b\in B$ and $f\in L^1(G)$. Using the easily verifiable identity
$$
W^{\omega*}(\rho^\omega_g\otimes1)W^\omega=\rho_g\otimes\lambda^\omega_g,
$$
we get the required isomorphism
$$
\overline{\alpha(B)(1\otimes\rho^\omega(C^*_r(G,\omega)))}\to A_\omega,\ \ x\mapsto W^{\omega*}_{23}(x\otimes1)W^\omega_{23}.
$$
In order to see that this isomorphism respects the coactions, we just have to check that
$$
\delta^\omega(1\otimes \rho_g\otimes \lambda^\omega_g)=1\otimes \rho_g\otimes \lambda^\omega_g\otimes\lambda_g,
$$
that is, $W(\lambda^\omega_g\otimes1)W^*=\lambda^\omega_g\otimes\lambda_g$. But this follows immediately from $W(\lambda_g\otimes1)W^*=\lambda_g\otimes\lambda_g$, since~$\lambda^\omega_g$ is $\lambda_g$ multiplied by a function that automatically commutes with the first leg of $W$.
\ep

Consider now an arbitrary coaction $\delta$ of $G$ on a C$^*$-algebra $A$ and choose two cocycles $\omega,\nu\in Z^2(G;\T)$. Using the coaction $\delta^\omega$ on $A_\omega$ we can define the $\nu$-deformation $(A_\omega)_\nu$ of $A_\omega$.

\begin{proposition} \label{piterate}
The map $$A_{\omega\nu}\to M(A\otimes K\otimes K),\ \ x\mapsto W_{23}\tilde\nu_{23}^*(x\otimes1)\tilde\nu_{23}W^*_{23},$$ defines an isomorphism $A_{\omega\nu}\cong (A_\omega)_\nu$.
In particular, the map $\eta^\omega\delta\colon A\to M(A\otimes K\otimes K)$ defines an isomorphism $A\cong (A_\omega)_{\bar\omega}$.
\end{proposition}

\bp For $a\in A$ and $\varphi\in K^*$ consider the element
$$
x=(\iota\otimes\iota\otimes\varphi)\eta^\omega\delta(a)
=(\iota\otimes\iota\otimes\varphi)(\iota\otimes(\lambda^\omega\boxtimes\lambda^{\bar\omega}))\delta(a)
\in A_\omega.
$$
Recall that $\lambda^\omega\boxtimes\lambda^{\bar\omega}$ denotes the representation of $C^*_r(G)$ defined by $\lambda_g\mapsto \lambda^\omega_g\otimes\lambda^{\bar\omega}_g$. Then
$$
\delta^\omega(x)=W_{23}(x\otimes1)W_{23}^*=(\iota\otimes\iota\otimes\varphi\otimes\iota)(W_{24}
((\iota\otimes(\lambda^\omega\boxtimes\lambda^{\bar\omega}))\delta(a)\otimes1)W_{24}^*).
$$
Since $W(\lambda^\omega_g\otimes1)W^*=\lambda^\omega_g\otimes\lambda_g$, as was already used in the proof of the previous proposition, the above identity can be written as
$$
\delta^\omega(x)=(\iota\otimes\iota\otimes\varphi\otimes\iota)
(\iota\otimes((\lambda^\omega\boxtimes\lambda^{\bar\omega})\boxtimes\lambda))\delta(a).
$$
It follows that
$$
\eta^\nu\delta^\omega(x)=(\iota\otimes\iota\otimes\varphi\otimes\iota\otimes\iota)
(\iota\otimes((\lambda^\omega\boxtimes\lambda^{\bar\omega})\boxtimes
(\lambda^\nu\boxtimes\lambda^{\bar\nu})))\delta(a).
$$
Therefore $(A_\omega)_\nu$ is the closed linear span of elements of the form
$$
(\iota\otimes\iota\otimes\varphi\otimes\iota\otimes\psi)
(\iota\otimes((\lambda^\omega\boxtimes\lambda^{\bar\omega})
\boxtimes(\lambda^\nu\boxtimes\lambda^{\bar\nu})))\delta(a),
$$
where $a\in A$ and $\varphi,\psi\in K^*$.

Observe next that
$$
W\tilde\nu^*(\lambda_g^{\omega\nu}\otimes1)\tilde\nu W^*=\lambda^\omega_g\otimes\lambda^\nu_g,
$$
which is simply identity \eqref{ebasic} for the cocycle $\bar\nu$ multiplied on the left by $\tilde\omega(g^{-1},\cdot)\tilde\nu(g^{-1},\cdot)\otimes1$. It follows that the unitary
$$
\Sigma_{23}(\tilde\nu W^*\otimes \tilde{\nu}^*W^*)\Sigma_{23}\ \ \text{on}\ \ L^2(G)^{\otimes 4},
$$
where $\Sigma$ is the flip, intertwines the representation $(\lambda^\omega\boxtimes\lambda^{\bar\omega})\boxtimes(\lambda^\nu\boxtimes\lambda^{\bar\nu})$ of $C^*_r(G)$ with the representation $(\lambda^{\omega\nu}\boxtimes\lambda^{\bar\omega\bar\nu})\otimes1\otimes1$. Furthermore, for any $y\in C^*_r(G)$ we have
$$
(\Ad \tilde\nu W^* )(\iota\otimes\varphi\otimes\iota\otimes\psi)
((\lambda^\omega\boxtimes\lambda^{\bar\omega})\boxtimes(\lambda^\nu\boxtimes\lambda^{\bar\nu}))(y)
=\varpi_{24}((\lambda^{\omega\nu}\boxtimes\lambda^{\bar\omega\bar\nu})(y)\otimes1\otimes1),
$$
where $\varpi=(\varphi\otimes\psi)(\Ad W\tilde{\nu})\in (K\otimes K)^*$. Therefore for any $a\in A$ we get
$$
(\Ad \tilde\nu_{23} W^*_{23})(\iota\otimes\iota\otimes\varphi\otimes\iota\otimes\psi)
(\iota\otimes((\lambda^\omega\boxtimes\lambda^{\bar\omega})
\boxtimes(\lambda^\nu\boxtimes\lambda^{\bar\nu})))\delta(a)
=\varpi_{35}(\eta^{\omega\nu}\delta(a)\otimes1\otimes1).
$$
This shows that $\Ad (\tilde\nu_{23} W^*_{23})$ maps the algebra $(A_\omega)_\nu$ onto $A_{\omega\nu}\otimes1$, which proves the first part of the proposition. Then the second part also follows, since the deformation of $A$ by the trivial cocycle is equal to $\delta(A)$.
\ep


\bigskip

\section{K-theory}

We say that two cocycles $\omega_0,\omega_1\in Z^2(G;\T)$ are homotopic if there exists a $C([0,1];\T)$-valued Borel $2$-cocycle $\Omega$ on $G$ such that $\omega_i=\Omega(\cdot,\cdot)(i)$ for $i=0,1$. Our goal is to show that under certain assumptions on $G$ the deformed algebras $A_{\omega_0}$ and $A_{\omega_1}$ have isomorphic $K$-theory. For this we will use the following slight generalization of invariance under homotopy of cocycles of $K$-theory of reduced twisted group C$^*$-algebras, proved in~\cite{ELPW}.

\begin{theorem} \label{tktheory}
Assume $G$ satisfies the Baum-Connes conjecture with coefficients. Then for any action $\alpha$ of $G$ on a C$^*$-algebra $B$ and any two homotopic cocycles $\omega_0,\omega_1\in Z^2(G;\T)$, for the corresponding reduced twisted crossed products we have $K_*(B\rtimes_{\alpha,\omega_0}G)\cong K_*(B\rtimes_{\alpha,\omega_1}G)$.
\end{theorem}

The proof follows the same lines as that of \cite[Theorem~1.9]{ELPW}. The starting point is the isomorphism
$$
K\otimes(B\rtimes_{\alpha,\omega}G)\cong(K\otimes B)\rtimes_{\Ad\rho^{\bar\omega}\otimes\alpha}G,\ \ x\mapsto \omega_{13}^*V_{13}x V_{13}^*\omega_{13},
$$
which maps $\rho^{\bar\omega}_g\otimes1\otimes\rho^\omega_g$ into $1\otimes1\otimes\rho_g$. This is a particular case of the Packer-Raeburn stabilization trick, see \cite[Section~3]{PR}. Therefore instead of twisted crossed products we can consider $(K\otimes B)\rtimes_{\Ad\rho^{\bar\omega}\otimes\alpha}G$.

Now, given a homotopy $\Omega$ of cocycles, consider the action $\Ad\rho^{\bar\Omega}$ of $G$ on $C[0,1]\otimes K$ defined, upon identifying $C[0,1]\otimes K$ with $C([0,1];K)$, by $(\Ad\rho^{\bar\Omega}_g)(f)(t)=(\Ad\rho^{\bar\omega_t}_g)(f(t))$, where $\omega_t=\Omega(\cdot,\cdot)(t)$.

\begin{lemma}[cf.~Proposition~1.5 in \cite{ELPW}] \label{lexterior}
For any compact subgroup $H\subset G$ and any $t\in[0,1]$, the restrictions of the actions $\Ad\rho^{\bar\Omega}$ and $\operatorname{id}\otimes\Ad\rho^{\bar\omega_t}$ to $H$ are exterior equivalent.
\end{lemma}

Note that this is easy to see for homotopies of the form $\omega_t=\omega_0e^{itc}$ usually considered in applications, where $c$ is an $\R$-valued Borel $2$-cocycle. Indeed, by \cite[Theorem~2.3]{Moore1} the second cohomology of a compact group with coefficients in $\R$ is trivial, so there exists a Borel function $b\colon H\to\R$ such that $c(h',h)=b(h')+b(h)-b(h'h)$. Extend $b$ to a function on $G$ as follows. Choose a Borel section $s\colon G/H\to G$ of the quotient map $G\to G/H$, $g\mapsto\dot{g}$, such that $s(\dot{e})=e$. Then put
$$
b(g)=b(s(\dot{g})^{-1}g)-c(s(\dot{g}),s(\dot{g})^{-1}g)+b(e).
$$
A simple computation shows that $c(g,h)=b(g)+b(h)-b(gh)$ for all $g\in G$ and $h\in H$. Then the unitaries $u_h\in M(C[0,1]\otimes K)$ defined by $u_h(t)=e^{it(b-b(\cdot h))}$ form a $1$-cocycle for the action $(\Ad\rho^{\bar\Omega})|_H$ such that $\Ad (u_h\rho^{\bar\Omega}_h)=\operatorname{id}\otimes\Ad\rho^{\bar\omega_0}_h$.

\bp[Proof of Theorem~\ref{tktheory}] For every $t\in[0,1]$ consider the evaluation map $\ev_t\colon C[0,1]\otimes K\otimes B\to K\otimes B$. It is $G$-equivariant with respect to the action $\Ad\rho^{\bar\Omega}\otimes\alpha$ of $G$ on $C[0,1]\otimes K\otimes B$ and the action $\Ad\rho^{\bar\omega_t}\otimes\alpha$ of $G$ on $K\otimes B$. We claim that it induces an isomorphism
$$
(\ev_t\rtimes G)_*\colon K_*((C[0,1]\otimes K\otimes B)\rtimes_{\Ad\rho^{\bar\Omega}\otimes\alpha}G)\to K_*((K\otimes B)\rtimes_{\Ad\rho^{\bar\omega_t}\otimes\alpha}G).
$$
By \cite[Proposition~1.6]{ELPW} in order to show this it suffices to check that for every compact subgroup $H$ of $G$ the map $\operatorname{ev}_t$ induces an isomorphism
$$
(\ev_t\rtimes H)_*\colon K_*((C[0,1]\otimes K\otimes B)\rtimes_{\Ad\rho^{\bar\Omega}\otimes\alpha}H)\to K_*((K\otimes B)\rtimes_{\Ad\rho^{\bar\omega_t}\otimes\alpha}H).
$$
By Lemma~\ref{lexterior} the action $\Ad\rho^{\bar\Omega}\otimes\alpha$ of $H$ on $C[0,1]\otimes K\otimes B$ is exterior equivalent to the action $\operatorname{id}\otimes\Ad\rho^{\bar\omega_t}\otimes\alpha$, so that
$$
(C[0,1]\otimes K\otimes B)\rtimes_{\Ad\rho^{\bar\Omega}\otimes\alpha}H\cong C[0,1]\otimes((K\otimes B)\rtimes_{\Ad\rho^{\bar\omega_t}\otimes\alpha}H).
$$
If the cocycle $U=\{u_h\}_{h\in H}$ defining the exterior equivalence is chosen such that $u_h(t)=1$ for all~$h\in H$, then the corresponding homomorphism
$$
C[0,1]\otimes((K\otimes B)\rtimes_{\Ad\rho^{\bar\omega_t}\otimes\alpha}H)\to (K\otimes B)\rtimes_{\Ad\rho^{\bar\omega_t}\otimes\alpha}H
$$
is simply the evaluation at $t$. Obviously, it defines an isomorphism in $K$-theory.
\ep

Combining Theorems~\ref{tmain2} and~\ref{tktheory} we get the following result that generalizes several earlier results in the literature~\cite{Ri2,Yam}.

\begin{corollary}
Assume $G$ satisfies the Baum-Connes conjecture with coefficients. Then for any coaction $\delta$ of $G$ on a C$^*$-algebra $A$ and any two homotopic cocycles $\omega_0,\omega_1\in Z^2(G;\T)$, we have an isomorphism $K_*(A_{\omega_0})\cong K_*(A_{\omega_1})$.
\end{corollary}

We finish by noting  that for some groups it is possible to prove a stronger result. For example, generalizing Rieffel's result for $\R^d$~\cite{Ri2} we have the following.

\begin{proposition}
If $G$ is a simply connected solvable Lie group, then for any coaction $\delta$ of $G$ on a C$^*$-algebra $A$ and any cocycle $\omega\in Z^2(G;\T)$ we have $K_*(A_\omega)\cong K_*(A)$.
\end{proposition}

\bp
By the stabilization trick and Connes' Thom isomorphism we have $K_i(A\rtimes_\delta\hat G\rtimes_{\hat\delta,\omega}G)\cong K_{i+\dim G}(A\rtimes_\delta \hat G)\cong K_i(A\rtimes_\delta\hat G\rtimes_{\hat\delta}G)\cong K_i(A)$.
\ep


\bigskip

\appendix

\section{Rieffel's deformation}

It was stated by Kasprzak~\cite{Kas} that for $G=\R^d$ his approach to deformation, which our construction extends, is equivalent to that of Rieffel~\cite{Ri1}, but no proof of this was given. A sketch of a possible proof was then proposed by Hannabuss and Mathai~\cite{HM}, but in our opinion it is not easy to obtain a complete proof following the suggested strategy. The goal of this appendix is to give a different rigorous proof using completely positive maps constructed by Kaschek, Neumaier and Waldmann~\cite{KNW}.

\smallskip

We will use the conventions in~\cite{KNW} that are slightly different from those of Rieffel. Assume $V$ is a $2n$-dimensional Euclidean space with scalar product $\langle\cdot,\cdot\rangle$, and $J$ is a complex structure on $V$, so $J$ is an orthogonal transformation and $J^2=-1$. Fix a deformation parameter $h>0$.

Assume we are given an action $\alpha$ of $V$ on a C$^*$-algebra $A$. Denote by $A_\infty$ the subalgebra of smooth vectors for this action. It is a Fr\'echet algebra equipped with differential norms $\|\cdot\|_k$, $k\ge1$. Rieffel defines a new product $*_h$ on $A_\infty$ by
$$
a*_hb=\frac{1}{(\pi h)^{2n}}\int_{V\times V}\alpha_x(a)\alpha_y(b)e^{-\frac{2i}{h}\langle x,Jy\rangle}dx\,dy,
$$
where the integral is understood as an oscillatory integral.

Denote by $A_x$ the spectral subspace of~$A_\infty$ corresponding to $x\in V$, so $A_x$ consists of elements $a\in A$ such that $\alpha_y(a)=e^{i\langle x,y\rangle}a$ for all $y\in V$. Then for $a\in A_x$ and $b\in A_y$ we have
$$
a*_h b=e^{\frac{ih}{2}\langle x,Jy\rangle}ab.
$$
Note that the spectral subspaces are often trivial, so this formula by no means determines $*_h$. Nevertheless it indicates that the cocycle of deformation is $\omega(x,y)=e^{\frac{ih}{2}\langle x,Jy\rangle}$. The Rieffel deformation of~$A$ is a certain C$^*$-algebraic completion of $A_\infty$ equipped with the product $*_h$ and with the involution inherited from~$A$, see~\cite{Ri1} for details. We denote it by $\tilde A_\omega$.

The action $\alpha$ can be viewed as a coaction $\delta$ of $V$ on $A$. Namely, define the Fourier transform $\F\colon L^2(V)\to L^2(V)$ by
$$
(\F f)(x)=\frac{1}{(2\pi)^n}\int_Vf(y)e^{-i\langle x,y\rangle}dy.
$$
Then $\Ad\F$ defines an isomorphism of $C_0(V)$ onto $C^*_r(V)$, and by letting $\delta=\Ad(1\otimes\F)\alpha$ we get a coaction of $V$ on $A$. Note that $a\in A$ lies in the spectral subspace $A_x$ if and only if $\delta(a)=a\otimes\lambda_x$, in agreement with our previous notation. We can then consider the $\omega$-deformation $A_\omega$ of $A$. Our aim is to construct an isomorphism between~$A_\omega$ and $\tilde A_\omega$.

\smallskip

Following~\cite{KNW} define a map $\Phi\colon A\to A$ by
$$
\Phi(a)=\frac{1}{(\pi h)^n}\int_Ve^{-\frac{1}{h}\|x\|^2}\alpha_x(a)dx.
$$
We have
\begin{equation} \label{ephi1}
\Phi(a)=e^{-\frac{h}{4}\|x\|^2}a\ \ \text{for}\ \ a\in A_x.
\end{equation}
The image of $\Phi$ is contained in $A_\infty$. So we can consider $\Phi$ as a map $\tilde T\colon A\to\tilde A_\omega$. Identifying $A$ with Rieffel's deformation of $\tilde A_\omega$ corresponding to the complex structure $-J$, we also get a similarly defined map $\tilde S\colon\tilde A_\omega\to A$, so the restriction of $\tilde S$ to $A_\infty$ coincides with the restriction of $\Phi$ to $A_\infty$. Since $\Phi$ considered as a map $(A,\|\cdot\|)\to(A_\infty,\|\cdot\|_k)$ is bounded for any $k$, the map $\tilde T\colon A\to \tilde A_\omega$ is bounded by standard estimates for the operator norm on $\tilde A_\omega$, see \cite[Proposition~4.10]{Ri1}. By symmetry the map $\tilde S$ is also bounded. The main result in \cite{KNW} states that the maps $\tilde T$ and $\tilde S$ are completely positive. We will reprove this a bit later.

We want to define analogues of the maps $\tilde T$ and $\tilde S$ for $A_\omega$. For this, define a unit vector $\xi_0\in L^2(V)$~by
$$
\xi_0(x)=\left(\frac{h}{2\pi}\right)^{n/2}e^{-\frac{h}{4}\|x\|^2}.
$$
Consider the normal state $\varphi_0=(\cdot\,\xi_0,\xi_0)$ on $B(L^2(V))$. We have
\begin{equation*}
\varphi_0(\lambda^\omega_x)=\varphi_0(\lambda^{\bar\omega}_x)=e^{-\frac{h}{4}\|x\|^2}.
\end{equation*}
Note that this means that on the C$^*$-algebra generated by the operators~$\lambda^\omega_x$, which is the algebra of canonical commutation relations for the space $V$ equipped with the Hermitian scalar product $h\langle x,y\rangle+ih\langle x,Jy\rangle$, the state $\varphi_0$ is simply the vacuum state.

Define $T\colon A\to A_\omega$ and $S\colon A_\omega\to A$ by
$$
T(a)=(\iota\otimes\iota\otimes\varphi_0)\eta^\omega\delta(a),\ \ S(b)=(\iota\otimes\varphi_0)(b).
$$
Using that $\delta(A)(1\otimes C^*_r(V))\subset A\otimes C^*_r(V)$ it is not difficult to see that the image of $S$ is indeed contained in $A$ rather than in $M(A)$. This will also become clear from the proof of Lemma~\ref{la2} below.

The maps $T$ and $S$ are completely positive. Using that $\eta^\omega\delta(a)=a\otimes\lambda^\omega_x\otimes\lambda^{\bar\omega}_x$ for $a\in A_x$, we get
\begin{equation} \label{ephi2}
T(a)=e^{-\frac{h}{4}\|x\|^2}a\otimes\lambda^\omega_x\ \ \text{and}\ \ S(a\otimes\lambda^\omega_x)=e^{-\frac{h}{4}\|x\|^2}a\ \ \text{for}\ \ a\in A_x.
\end{equation}

\begin{lemma} \label{la1}
For any $n\ge1$ and $a_1,\dots,a_n\in A$ we have
$$
\tilde S(\tilde T(a_1)\dots \tilde T(a_n))=S(T(a_1)\dots T(a_n)).
$$
\end{lemma}

\bp If for every $j$ the element $a_j$ lies in a spectral subspace $A_{x_j}$, then the identity in the formulation follows immediately from \eqref{ephi1} and \eqref{ephi2}. We will show that this is enough to conclude that it holds for arbitrary elements.

\smallskip

We claim that there exists a von Neumann algebra $M$ containing $A$ such that the action $\alpha$ of $V$ on $A$ extends to a continuous (in the von Neumann algebraic sense) action of $V$ on $M$ and such that~$M$ is generated as a von Neumann algebra by the spectral subspaces of this action. Indeed, first represent the crossed product $A\rtimes_\alpha V$ faithfully on some Hilbert space $H$ and consider the von Neumann algebra $N\subset B(H)$ generated by $A$. The action $\alpha$ of $V$ on $A$ extends to an action $\beta$ of~$V$ on~$N$. Consider the double crossed product $M=N\rtimes_\beta V\rtimes_{\hat\beta} \hat V$ in the von Neumann algebraic sense. By the Takesaki duality we have $(M,\hat{\hat\beta})\cong (N\bar\otimes B(L^2(V)),\beta\otimes\Ad\rho)$. This gives us an equivariant embedding of $A\subset N$ into $M$ equipped with the action $\hat{\hat\beta}$. It is also clear that $M$ is generated by the spectral subspaces of the action, so our claim is proved.

\smallskip

We continue to denote by $\alpha$ the action of $V$ on $M$. Denote by $\M\subset M$ the set of elements $a\in M$ such that the map $x\mapsto\alpha_x(a)$ is norm-continuous. This is an ultrastrongly operator dense C$^*$-subalgebra of $M$. We continue to denote by $T,S,\tilde T,\tilde S$ the maps defined for the C$^*$-algebra $\M$ in place of $A$. The maps $T$ and $S$ have obvious extensions to normal maps between the von Neumann algebras generated by $\M$ and $\M_\omega$. On the other hand, the map $\Phi$,
$$
\Phi(a)=\frac{1}{(\pi h)^n}\int_Ve^{-\frac{1}{h}\|x\|^2}\alpha_x(a)dx,
$$
is still well-defined on $M$, but now the integral should be taken with respect to the ultrastrong operator topology. The image of $M$ under $\Phi$ is contained in $\M_\infty$. It therefore makes sense to ask whether the identity
$$
\Phi(\Phi(a_1)*_h\dots *_h\Phi(a_n))=S(T(a_1)\dots T(a_n))
$$
holds for all $a_1,\dots,a_n\in M$, which would imply the assertion of the lemma. Since this identity holds for $a_1,\dots,a_n$ lying in spectral subspaces of $M$, it suffices to show that both sides of the identity are normal maps in every variable $a_j$ running through the unit ball $M_1$ of $M$. This is clearly the case for the right hand side. In order to prove the same for the left hand side it suffices to show that for any $b,c\in\M_\infty$ the map
$$
M_1\to M,\ \ a\mapsto \Phi(b*_h\Phi(a)*_hc),
$$
is continuous in the ultrastrong operator topology.

Using basic estimates for oscillatory integrals, see \cite[Chapter~1]{Ri1}, and the fact that the map $\Phi$ is bounded as a map $(M,\|\cdot\|)\to (\M_\infty,\|\cdot\|_k)$ for every $k$, it is easy to check that $\Phi(b*_h\Phi(a)*_hc)$ can be approximated in norm uniformly in $a\in M_1$ by integrals of the form
$$
\int_{V^3}\psi(x,y,z)\alpha_x(b)\alpha_y(a)\alpha_z(c)dx\,dy\,dz,
$$
where $\psi$ is a smooth compactly supported function and the integral is taken with respect to the ultrastrong operator topology. Since such integrals are clearly continuous in $a\in M_1$ in this topology, this finishes the proof of the lemma.
\ep

We will need the above lemma only for $n=1,2,3$.

\begin{lemma} \label{la2}
The maps $\tilde T$ and $\tilde S$ are completely positive, and all four maps $T,S,\tilde T,\tilde S$ are injective and their images are dense.
\end{lemma}

\bp We begin by proving that $\tilde T$ and $\tilde S$ are injective. It suffices to consider $\tilde T$. Assume $\tilde T(a)=0$. Then
\begin{equation} \label{ea1}
\frac{1}{(\pi h)^n}\int_Ve^{-\frac{1}{h}\langle x-y,x-y\rangle}\alpha_x(a)dx=\alpha_y(\Phi(a))=0
\end{equation}
for all $y\in V$, hence, by analyticity, for all $y$ in the complexification $V_\C$ of $V$. This implies that the Fourier transform of the $A$-valued function $x\mapsto e^{-\frac{1}{h}\|x\|^2}\alpha_x(a)$ is zero, whence $a=0$.

\smallskip

Next we will show that the images of $\tilde T$ and $\tilde S$ are dense. It suffices to consider $\tilde S$, and then it is enough to show that the image of $\Phi$ is dense. It is well-known, and is easy to check using e.g.\ Wiener's Tauberian theorem, that the translations of the function $e^{-\frac{1}{h}\|x\|^2}$ span a dense subspace of~$L^1(V)$. Using the first equality in \eqref{ea1} and that $\alpha_y(\Phi(a))=\Phi(\alpha_y(a))$, we conclude that the closure of the image of~$\Phi$ contains all elements of the form $\int_Vf(x)\alpha_x(a)dx$ with $f\in L^1(V)$. Hence this closure coincides with~$A$.

\smallskip

Let us show now that $\tilde T$ and $\tilde S$ are completely positive. Again, it is enough to consider $\tilde S$. Since by Lemma~\ref{la1} we have
$$
\tilde S(\tilde T(a)^*\tilde T(a))=S(T(a)^*T(a))\ge0,
$$
and the image of $\tilde T$ is dense, we see that $\tilde S$ is positive. Passing to deformations of matrix algebras over $A$ we conclude that $\tilde S$ is completely positive. This finishes the proof of the lemma for $\tilde T$ and $\tilde S$.

\smallskip

Turning to $T$ and $S$, by Lemma~\ref{la1} we have $ST=\tilde S\tilde T=\Phi^2$. Since the map $\Phi$ is injective and its image is dense, it follows that the map $T$ is injective and the image of~$S$ is dense. Consider the maps $T'\colon A_\omega\to (A_\omega)_{\bar\omega}$ and $S'\colon (A_\omega)_{\bar\omega}\to A_\omega$ defined by $(A_\omega,\delta^\omega)$ in the same way as $T$ and $S$ were defined by $(A,\delta)$. Then $T'$ is injective and the image of $S'$ is dense. By Proposition~\ref{piterate} the map~$\eta^\omega\delta$ defines an isomorphism $A\cong (A_\omega)_{\bar\omega}$. By definition of $T$ and $S'$ we immediately get $T=S'\eta^\omega\delta$. Hence the image of $T$ is dense. We also have $\eta^\omega\delta S=T'$. Indeed, a simple computation similar to the ones used in the proof of Proposition~\ref{piterate} shows that for $b=(\iota\otimes\iota\otimes\varphi)\eta^\omega\delta(a)\in A_\omega$ we have
$$
\eta^\omega\delta S(b)=(\iota\otimes\iota\otimes\iota\otimes\varphi_0\otimes\varphi)
(\iota\otimes((\lambda^\omega\boxtimes\lambda^{\bar\omega})
\boxtimes(\lambda^\omega\boxtimes\lambda^{\bar\omega})))\delta(a)
$$
and
$$
T'(b)=(\iota\otimes\iota\otimes\varphi\otimes\iota\otimes\varphi_0)
(\iota\otimes((\lambda^\omega\boxtimes\lambda^{\bar\omega})
\boxtimes(\lambda^{\bar\omega}\boxtimes\lambda^\omega)))\delta(a).
$$
Alternatively, the identity $\eta^\omega\delta S=T'$ is immediate on elements of the form $a\otimes\lambda_x$, where $a\in A_x$, hence it holds on arbitrary elements by an argument similar to the one used in the proof of Lemma~\ref{la1}. It follows that the map $S$ is injective.
\ep

Note that instead of injectivity we will only need to know that $S$ and $\tilde S$ are faithful. While it is obvious that $\tilde S$ is faithful, this is not the case for $S$, since the state $\varphi_0$ is very far from being faithful on the von Neumann algebra generated by $C^*_r(V,\omega)$. This von Neumann algebra is a factor of type~I$_\infty$ and~$\varphi_0$ is a normal pure state on it, as can be shown by recalling that $\varphi_0$ defines the vacuum state on the algebra of canonical commutation relations generated by the operators $\lambda^\omega_x$.

\begin{theorem}
There is a unique isomorphism $A_\omega\cong\tilde A_\omega$ that maps $T(a)$ into $\tilde T(a)$ for every $a\in A$.
\end{theorem}

\bp Assume first that there exists a faithful state $\psi$ on $A$. Consider the positive linear functionals $\psi_\omega=\psi S$ and $\tilde\psi_\omega=\psi \tilde S$ on $A_\omega$ and $\tilde A_\omega$. Since the positive maps $S$ and $\tilde S$ are faithful, these functionals are faithful. Consider the faithful GNS-representation of $A_\omega$ on $H$ with cyclic vector $\xi$ defining $\psi_\omega$, and the faithful GNS-representation of $\tilde A_\omega$ on $\tilde H$ with cyclic vector $\tilde\xi$ defining $\tilde\psi_\omega$. By Lemma~\ref{la1} for $n=2$ we have
$$
(T(a)\xi,T(b)\xi)=(\tilde T(a)\tilde\xi,\tilde T(b)\tilde\xi).
$$
Since the images of $T$ and $\tilde T$ are dense, it follows that there exists a unitary operator $U\colon H\to \tilde H$ such that $UT(a)\xi=\tilde T(a)\tilde\xi$. By Lemma~\ref{la1} for $n=3$ we have
$$
(T(a)T(b)\xi,T(c)\xi)=(\tilde T(a)\tilde T(b)\tilde\xi,\tilde T(c)\tilde\xi),
$$
that is, $(UT(a)T(b)\xi,UT(c)\xi)=(\tilde T(a)UT(b)\xi,UT(c)\xi)$. Therefore $UT(a)=\tilde T(a)U$, so $\Ad U$ defines the required isomorphism.

In the general case the proof is basically the same, but instead of one state $\psi$ we have to choose a faithful family of states on $A$ and consider direct sums of the GNS-representations defined by the corresponding positive linear functionals on $A_\omega$ and $\tilde A_\omega$.
\ep


\bigskip

\end{document}